\documentclass[12pt]{article}

\usepackage{amsmath,amsthm,amsfonts,amssymb}

\def\qed{{\unskip\nobreak\hfil\penalty50
\hskip2em\hbox{}\nobreak\hfil$\square$
\parfillskip=0pt \finalhyphendemerits=0\par}\medskip}
\def\proof{\trivlist \item[\hskip \labelsep{\bf Proof\ }]}
\def\endproof{\null\hfill\qed\endtrivlist\noindent}

\def\a{\alpha}
\def\b{\beta}
\def\de{\delta}
\def\e{\varepsilon}
\def\g{\gamma}

\def\l{\lambda}

\def\r{\rho}
\def\phi{\varphi}

\def\om{\omega}
\def\Om{\Omega}

\newtheorem{theorem}{Theorem}
\newtheorem{lemma}[theorem]{Lemma}
\newtheorem{corollary}[theorem]{Corollary}

\newtheorem{proposition}[theorem]{Proposition}

\def\setminus{\smallsetminus}
\def\A{{\cal A}}
\def\B{{\cal B}}
\def\C{{\cal C}}

\def\M{{\cal M}}
\def\N{{\cal N}}
\def\R{{\cal R}}
\def\L{{\cal L}}
\def\I{{\cal I}}
\def\H{{\cal H}}
\def\K{{\cal K}}
\def\S{{\cal S}}
\def\U{{\cal U}}

\def\emptyset{\varnothing}
\def\res{{\upharpoonright}}
\def\PSL{{{\rm PSL}(2,\mathbb R)}}
\def\tPSL{{{\rm PSL}(2,\mathbb R)}^{\tilde{}
}}

\title{{\bf
Conformal Subnets \\
and Intermediate Subfactors}}

\author{
{\sc Roberto Longo}\footnote{Supported in part by MURST and 
INDAM-GNAMPA.} \\ Dipartimento di Matematica\\ Universit\`a di Roma 
``Tor Vergata''\\ Via della Ricerca Scientifica, I-00133 Roma, Italy\\ 
e-mail: {\tt longo@mat.uniroma2.it}}
\begin{document}
\date{}
\maketitle
\bigskip
\centerline{\sl Dedicated to Rudolf Haag on the occasion of his
eightieth birthday}
\bigskip
\bigskip
\bigskip
\bigskip
\begin{abstract}
Given an irreducible local conformal net $\A$ of von Neumann algebras 
on $S^1$ and a finite-index conformal subnet $\B\subset\A$, we show 
that $\A$ is completely rational iff $\B$ is completely rational.  In 
particular this extends a result of F.~Xu for the orbifold 
construction.  By applying previous results of Xu, many coset models 
turn out to be completely rational and the structure results in 
\cite{KLM} hold.  Our proofs are based on an analysis of the net 
inclusion $\B\subset\A$; among other things we show that, for a fixed 
interval $I$, every von Neumann algebra $\R$ intermediate between 
$\B(I)$ and $\A(I)$ comes from an intermediate conformal net $\L$ 
between $\B$ and $\A$ with $\L(I)=\R$.  We make use of a theorem of 
Watatani (type II case) and Teruya and Watatani (type III case) 
on the finiteness of the set $\mathfrak I(\N,\M)$ 
of intermediate subfactors in an irreducible inclusion of factors 
$\N\subset\M$ with finite Jones index $[\M:\N]$. We provide 
a unified proof of this result that gives in particular an explicit 
bound for the cardinality of $\mathfrak I(\N,\M)$ which depends 
only on $[\M:\N]$.
\end{abstract}
\newpage
\section{Introduction}
Operator algebraic methods have been used to good effect in Conformal 
Quantum Field Theory, in particular in understanding general model 
independent structure (e.g. \cite{BGL,KLM,F,GL2,GLW,R}), in the 
analysis of concrete models (e.g. \cite{BE,X1,X3,X}) and for 
applications in different contexts (e.g. \cite{L4}).  In most cases it 
seems to be impossible to proceed by different methods.

Because of their relevance in different areas, among others Topological 
QFT and 3-manifold invariants, conformal models with
a rational and modular representation theory have been the
subject of much attention, also in the physical literature (cf. 
\cite{EK}).

In \cite{KLM} intrinsic, model independent conditions selecting a 
class of (local, irreducible) conformal nets $\A$ of von Neumann 
algebras on $S^1$ with the right rationality/modularity properties 
were given. $\A$ is \emph{completely rational} if
\begin{enumerate}
\item $\A$ is split,
\item $\A$ is strongly additive,
\item the 2-interval inclusion of factors $\A(E)\subset\A(E')'$ has 
finite Jones index $\mu_\A$.
\end{enumerate}
Here both $E\subset S^1$ and $E'\equiv S^1\setminus E$ are the union of 
two proper intervals.  The split and strongly additivity properties 
are well-studied basic properties, see Section \ref{Hered} for their 
definitions, and we do not dwell on them here, cf.
\cite{DL,BD,DLR} and \cite{BS,GLW}.
If $\A$ is completely rational, then 
$\A(E)\subset\A(E')'$ is obtained by a quantum double construction in 
\cite{LR}, in particular
\[
\mu_\A=\sum_i d(\r_i)^2\ ,
\]
where the sum is taken over all the irreducible sectors of $\A$.  Every 
representation of $\A$ (on a separable Hilbert space) is M\"{o}bius 
covariant and decomposes into the direct sum of irreducible 
representations with finite statistical dimension.  There are only 
finitely many inequivalent irreducible representations, i.e.  $\A$ is 
rational, and the associated braiding is non-degenerate, i.e.  the 
representation  tensor category is modular.

At this point the problem of verifying the complete rationality of 
known models arises.  Certain examples were discussed in \cite{KLM}.  
As an illustration from \cite{KLM}, consider the case of a non-trivial 
finite group $G$ acting on a completely rational $\A$; if the 
fixed-point {\it orbifold} subnet $\A^G$ is also completely rational, 
then $\mu_{\A^G}\geq |G|^2$, while $\sum_{\pi\in\hat G}d(\r_{\pi})^2=|G|$, 
where the the $\r_{\pi}$'s are the untwisted DHR sectors of $\A^G$ 
\cite{DHR}, and this shows that twisted sectors must appear.  As $\A$ 
is the initial data, one would infer the complete rationality of 
$\A^G$ from that of $\A$.  By \cite{KLM} $\A^G$ inherits from $\A$ the 
split property and the finiteness of the $\mu$-index.  F. Xu \cite{X} 
has then shown that $\A^G$ also inherits the strong additivity 
property and this has inspired our paper.

We shall now show that if $\B$ is any
conformal subnet of $\A$ with finite index, then $\B$ is completely 
rational iff $\A$ is completely rational.

As a consequence, if $\B$ is a cofinite subnet of $\A$, 
namely $[\A:\B\vee\B^c]<\infty$, where $\B\vee\B^c$ is the subnet 
generated by $\B$ and its relative commutant $\B^c$ in $\A$, then $\A$ is 
completely rational iff both $\B$ and $\B^c$ are completely rational.

The subnet $\B^c$ is called the {\it coset} subnet associated with 
$\B\subset\A$, as it generalizes a coset construction that 
plays an important r\^{o}le in the theory of 
Kac-Moody Lie algebras, allowing one to construct of 
the minimal series representations of the Virasoro algebra 
\cite{GKO}. 

Coset models have been intensively studied by Xu in \cite{X1,X3} by 
operator algebraic methods. In one approach he makes use of \cite{KLM} 
too.  Thanks to his work, coset models associated 
with many loop group inclusions are cofinite, rational and modular, 
see the list in Section \ref{list}. Property 3 holds, but the 
validity of strong additivity was left open.

By our work in all these examples $\B^c$ is strongly additive, thus
$\B^c$ turns out to be completely rational and this 
completes the above discussion and explains the 
rationality/modularity structure better.

We now comment on our proof that the complete 
rationality property (and also the `split \& strongly additivity' 
property) for finite-index inclusions of conformal nets $\B\subset\A$
are hereditary.  
That `Property 1 \& 3' passes from $\A$ to $\B$ and viceversa is shown 
in \cite{KLM}.  The remaining more difficult point we have to prove is 
that $\B$ is strongly additive if $\A$ is split and strongly additive, 
see Sect.  \ref{Hered}.

To this end, we have analyzed a finite-index inclusion of conformal 
nets $\B\subset\A$ by considering the relative superselection 
structure.  In particular we show that, for a fixed interval $I$, 
every von Neumann algebra $\R$ intermediate between $\B(I)$ and 
$\A(I)$ comes from an intermediate conformal net $\L$ between $\A$ and 
$\B$ with $\L(I)=\R$.

Here we make use of a result of Watatani \cite{W} (in the type II case), 
following previous 
work by Popa \cite{P}, and Teruya and Watatani \cite{TW} (in the type 
III case) to the effect that the set $\mathfrak I(\N,\M)$ 
of intermediate subfactors in an irreducible inclusion of 
factors $\N\subset\M$ with finite Jones index $[\M:\N]$ is 
finite. We give a direct general proof of this result that works for 
arbitrary factors. This proof provides for the first time an explicit bound 
for the cardinality $|\mathfrak I(\N,\M)|$ of $\mathfrak I(\N,\M)$ which 
depends only on $[\M:\N]$ and implies that
\begin{equation*}
|\mathfrak I(\N,\M)|\leq \ell^{\ell} \ ,
\end{equation*}
where $\ell = [\M:\N]^2$.  There may be better bounds taking account 
of further structure associated with intermediate subfactors 
\cite{B,BJ}, see the comments at the end of Section 
\ref{subsectionbound}.

We conclude this introduction and include references to the books \cite{A,EK,H,T}
for basic facts on Operator Algebras and Quantum Field Theory, see
also \cite{K} for subfactors and sectors.
\section{On subfactors and intermediate factors}
The first part of this paper is devoted to an analysis of subfactors 
and intermediate factors, that will be used later on.
\subsection{Some basic structure}
Let $\N\subset\M$ be an irreducible inclusion of infinite factors with 
finite index $[\M:\N]$. We denote by $\g$ the canonical endomorphism 
of $\M$ into $\N$ and by $\theta$ the dual canonical endomorphism 
$\theta\equiv\g\res_\N$. The Q-system associated with $\g$ is denoted by 
$(\g,T,S)$, namely $T\in\M$ and $S\in\N$ are the unique (up to a 
phase) isometries in $(\iota,\g)$ and $(\iota,\theta)$, where 
$\iota$ always denotes the identity automorphism.

Let $\{[\rho_i]$, $i=0,\dots N\}$ be the family of the irreducible 
sectors in the decomposition of $[\theta]$, namely
\begin{equation}\label{decomposition}
[\theta] = \bigoplus_{i=0}^N N_i [\rho_i]\ .
\end{equation}
By Frobenious reciprocity for each $i$ the Hilbert space 
of isometries in $\M$ (not necessarily with right support $1$) 
\begin{equation}\label{Ki}
K_i\equiv\{R\in\M: Rx=\rho_i(x)R\ \forall x\in\N\}
\end{equation}
has dimension $N_i$, indeed the map
\begin{equation}\label{frob}
v\in (\r,\theta)\to v^* T \in K_i
\end{equation}
is an anti-linear isomorphism of $(\r,\theta)$ with $K_i$, whose inverse 
is given by
\[
X\in K_i\to [\M:\N]\e(TX^*)\in (\r,\theta)\ ,
\]
where $\e$ is the expectation of $\M$ onto $\N$. See \cite{LR,ILP} 
for the following.
\begin{lemma}\label{lemmaexpansion}
Let $\{R_{i,k}\}_{k=1}^{N_i}$ be an orthogonal basis of $K_i$ 
with the normalization $R_{i,k}^*R_{i,k}=d(\rho_i)$. 
Then every $X\in \M$ has a unique Fourier expansion
\begin{equation}
\label{expansion}
X=\sum^N_{i=0}\sum_{k=1}^{N_i} x_{i,k}R_{i,k}
\end{equation}
where the coefficients $x_{i,k}$ belong to $\N$, indeed
$x_{i,k}=\e(XR_{i,k}^*)$.
\end{lemma}
\proof It is immediate that $\e(R_{i,k}R^*_{j,h})\in(\r_j,\r_i)$, thus 
$\e(R_{i,k}R^*_{j,h})=0$ if $i\neq j$. Let $\{v_1,\dots v_{N_i}\}$ be 
an orthonormal basis of isometries for $K_i$, then 
$R_{i,k}\equiv [\M:\N]^{\frac{1}{2}}v_k^*T$ satisfy
\[
R_{i,k}^*R_{i,h}=[\M:\N]T^*v_k v_h^*T =\de_{hk} d(\r_i)
\]
because $v_k v_h^*\in(\theta,\theta)$ and $TT^*\in\M$ is the Jones 
projection for implementing the expectation $\e_1 :\N\to\theta(M)$ so 
that $\e_1\res_{(\theta,\theta)}$ is the associated trace, see 
\cite{L}. Therefore
\begin{equation}\label{orth}
\e(R_{i,k}R^*_{i,h})=[\M:\N]\e(v_k^*TT^*v_h)=v_k^* v_h=\de_{hk}.
\end{equation}
Now $\M=\N T$, thus $\M$ is generated by $\M=\sum_{i,k} \N\! R_{i,k}$ 
because $\sum K_i$ has right support one. By the orthogonality 
relations (\ref{orth}) every $X\in\M$ has the expansion given by 
formula (\ref{expansion}). 
\endproof
Denote by $\mathfrak I(\N,\M)$ the set of intermediate von Neumann 
algebras between $\N$ and $\M$. Clearly if $\R\in\mathfrak I(\N,\M)$, 
then $\R$ is a factor and indeed both $\N\subset\R$ and $\R\subset\M$ 
are irreducible finite-index inclusions of factors.

Let $\R\in\mathfrak I(\N,\M)$ and set $K'_i\equiv K_i\cap\R$. Then 
$K'_i$ is a Hilbert subspace of $K_i$, so we may choose the
$R_{i,k}$ so that $\{R_{i,k}\}_{k=1}^{N'_i}$ is a 
basis for $K'_i$, where $N'_i={\rm dim}K'_i$. We may also re-order 
the $\r_i$'s so that $N'_i > 0$ iff $i\leq N'$ where $N'\leq N$ is 
an integer.
\begin{proposition}
With the above notations, $X\in\M$ belongs to $\R$ iff in the 
expansion (\ref{expansion}) $x_{i,k}=0$ for all $k>N'_i$, namely
\begin{equation}
\label{expansion'}
X=\sum^{N'}_{i=0}\sum_{k=1}^{N'_i} x_{i,k}R_{i,k}\ .
\end{equation}
In particular $\R$ is generated by $\N$ and the $K'_i$'s.
\end{proposition}
\proof
Recall that $[\M:\R]<\infty$, so there exists an expectation $\e_\R 
:\M\to\R$. From the definition (\ref{Ki}) $\e_\R(K_i)\subset K_i$, 
thus $\e_\R\res_{K_i}\in B(K_i)$ is a norm one projection. Clearly
$\e_\R(K_i)=K'_i$ thus if $X\in \M$ has the expansion 
(\ref{expansion}) we have
\[
\e_\R(X)=\sum^{N'}_{i=0}\sum_{k=1}^{N'_i} x_{i,k}R_{i,k} \ ,
\]
which implies the statement in the proposition.
\endproof
The following theorem and its corollary are due to Watatani and 
Teruya-Watatani 
\cite{W, TW}, related results are contained in \cite{P}. 
\begin{theorem}\label{Watatani} 
$\mathfrak I(\N,\M)$ is a finite set.
\end{theorem}
\begin{corollary}\label{invar}
Let $\L$ be an intermediate subfactor between $\N$ and $\M$ and 
$\b:\mathfrak G\to {\rm Aut}(\M)$ a (pointwise weakly continuous) 
action of a connected topological group $\mathfrak G$ with 
$\b_g(\N)=\N$, $g\in\mathfrak G$. Then 
$\b_g(\L)=\L$, $g\in\mathfrak G$.
\end{corollary}
\proof 
We consider on $\mathfrak I(\M,\N)$ the topology of pointwise weak 
convergence of the associated conditional expectations ($\L_i\to\L$ 
iff $\e_{\L_i}(x)\to\e_{\L}(x)$ weakly for all $x\in\M$).  Then $\b$ 
implements a continuous action of $\mathfrak G$ on $\mathfrak 
I(\M,\N)$.  The corollary is thus immediate because any continuous 
action of a connected group on a discrete set is trivial. 
\endproof
\noindent
A direct proof of Th. \ref{Watatani} for 
factors of arbitrary type will be given in next section, where we shall 
obtain in particular a bound for the cardinality $|\mathfrak I(\N,\M)|$ of 
$\mathfrak I(\N,\M)$. 
\subsection{A bound for the number of intermediate subfactors}
\label{subsectionbound}
Let $\N\subset\M$ be an irreducible inclusion with finite index and 
denote by $\e$ the conditional expectation from $\M$ to $\N$. 
We shall now determine a bound for $|\mathfrak I(\M,\N)|$.
Our proof is inspired by the papers \cite{C,W}.

We assume that there exists a faithful normal state $\omega$ on $\N$ 
(otherwise replacing it by a weight).  By considering the GNS 
representation of $\M$ associated with $\tilde\om\equiv\om\cdot\e$, we 
may assume that $\M$ acts on a Hilbert space $\H$ with cyclic and 
separating vector $\Om$ so that $\tilde\om=(\Omega,\cdot\Omega)$.  
Then $e\equiv [\N\Omega]$, equal to $JeJ$, is the Jones projection for 
$\N\subset\M$ and $\M_1\equiv\langle \M,e\rangle=J\N'J$ is the Jones 
extension, where $J$ is the modular conjugation of $\M$ associated 
with $\Om$.  The projection $e\in\N'\cap\M_1$ and
\[
\e(e)=\lambda\equiv[\M:\N]^{-1}.
\]
Let $\R$, $\S$ be intermediate factors between $\N$ and $\M$ and denote 
by $p\equiv [\R\Omega]$ and $q\equiv [\S\Omega]$ their associated 
Jones projections. Note that $JpJ =p$, $JqJ =q$ and 
$\R_1\equiv\langle \R,p\rangle=J\R'J$ and $\S_1\equiv\langle 
\S,q\rangle=J\S'J$ are the corresponding Jones extensions and so 
there is a chain of inclusions
\[
\N\subset \R,\S\subset\M\subset\R_1 ,\S_1\subset\M_1\ .
\]
Clearly $p\in\R'\cap\R_1$ and $q\in\S'\cap\S_1$, thus $p$ and $q$ both 
belong to $\N'\cap\M_1$.
\begin{proposition}
If $||p-q||<\lambda/2$ then $\R=\S$.
\end{proposition}
\proof
We may assume $\R\neq\N$ as otherwise $p=e$, thus $q=e$ because 
$q\geq e$ and $\l\leq 1$.

As we have $[\M:\N]=[\M:\R][\R:\N]$, it follows that
\begin{equation}\label{pp}
[\M:\N]\geq 2[\M:\R]
\end{equation}
because $[\R:\N]\geq 2$ \cite{J}. In particular $\lambda\leq 1/2$.
Let $\e_{\R_1}$ be the expectation from $\M_1$ onto $\R_1$ and set 
$q'\equiv\e_{\R_1}(q)$. Then $q'\in\S'\cap\R_1$ since obviously 
$q'\in\R_1$ and if $x\in\S$
\[
xq'=x\e_{\R_1}(q)=\e_{\R_1}(xq)=\e_{\R_1}(qx)=q'x
\]
because $\S\subset\R_1$. Moreover $0\leq q'\leq 1$ and $q'\neq 0$ 
because $\e_{\R_1}$ is positive and faithful. Setting $\de=\l/2$ we have
\begin{equation}\label{less}
||p-q'||=||p-\e_{\R_1}(q)||=||\e_{\R_1}(p-q)||\leq 
||p-q||<\de\ .
\end{equation}
Therefore the spectrum ${\rm sp}(q')\subset [0,\de)\cup (1-\de,1]$, see 
Lemma \ref{sp}.

Thus the spectral projection $q''\equiv \chi_{(1-\de,1]}(q')$ is a 
projection in $\S'\cap\R_1$ and $||q''-q'||<\de$, thus
\[
||p-q''||\leq ||p-q'||+||q'-q''||<2\de\leq 1\ ,
\]
and this implies that $p$ and $q''$ are equivalent projections of 
$\N'\cap\R_1$. Indeed the phase $v$ in the polar decomposition of 
$t\equiv pq''$ is a partial isometry in $\N'\cap\R_1$ with $v^* 
v=q''$, $vv^*=p$ (see Lemma \ref{polar} below). 
Then we can define an isomorphism $\Phi$ of $\S$ 
into $\R$ by
\[
\Phi(x)p\equiv vxv^*\ , \ x\in\S\ ,
\]
as $p\R_1 p=\R p$ \cite{J}. Moreover $\Phi(x)=x$ for all $x\in\N$ because 
$v\in\N'$. 

We have the intertwining relation
\[
\Phi(x)v = vx\ , \ x\in\S\ .
\]
With $\e'$ the conditional expectation from $\M_1$ onto $\M$ we then have
\[
\Phi(x)\e'(v) = \e'(v)x\ , \ x\in\S\ ,
\]
where $\e'(v)\in\N'\cap\M=\mathbb C$, thus $\e'(v)\neq 0$ would imply 
that $\Phi$ is the identity on $\S$ and $\S\subset\R$. Reversing the 
r\^{o}le of $\R$ and $\S$ also $\R\subset\S$, so $\R=\S$.

To show that indeed $\e'(v)\neq 0$ set $\lambda_0\equiv\e'(p)=[\M:\R]^{-1}$ 
and notice that, by using Lemma \ref{polar}, we have
\begin{multline*}
|\e'(v)-\lambda_0|\leq ||\e'(v)-\e'(t)||+
||\e'(t)-\e'(p)||\leq ||v-t||+||t-p||\\
=||v-t||+||p(q''-p)||
\leq 2\de +||q''-p||<4\de=2\l \ ,
\end{multline*}
thus $\e'(v)\neq 0$ because $2\l\leq\l_0$ by eq. (\ref{pp}).
\endproof
\begin{corollary}\label{corbound}
Let $\N\subset\M$ be an irreducible inclusion of factors. The 
cardinality of the set of intermediate factors between $\N$ and $\M$ 
is bounded by
\begin{equation}\label{bound}
|\mathfrak I(\N,\M)|\leq (4(n+2)\sqrt{n} + 1)^{n^2} \ ,	
\end{equation}
where $n$ is the largest integer such that $n+1\leq [\M:\N]$.
\end{corollary}
\proof 
By the above proposition $|\mathfrak I(\N,\M)|$ is dominated by the 
maximum number of projections $\geq e$ in $\N'\cap\M_1$ whose mutual 
distance is $\geq\l/2$.  As $e$ is a minimal central projection of 
$\N'\cap\M_1$, we can naturally embed $\N'\cap\M_1$ into $\mathbb 
C\oplus {\rm Mat}_m(\mathbb C)$, where $m$ an integer with $m+1\leq [\M:\N]$
($m=\sum^N_1 N_i$ in eq.  (\ref{decomposition})).  Indeed, as 
$J(\N'\cap\M_1)J=\N'\cap\M_1$, and Ad$J$ implements an 
anti-automorphism of $\N'\cap\M_1$, we can assume this 
anti-automorphism to extend to an anti-automorphism of $\mathbb C\oplus 
{\rm Mat}_m(\mathbb C)$ preserving the two components, in other words 
we may assume that the Ad$J$-invariant part of $\N'\cap\M_1$ is 
contained in $\mathbb R\oplus {\rm Mat}_m(\mathbb R)$.

Thus $|\mathfrak I(\N,\M)|$ is dominated by the maximum number of 
projections in ${\rm Mat}_m(\mathbb C)$ whose mutual distance is 
larger than $\frac{1}{2(n+2)}$ where $n$ is the largest interger
such that $n+1\leq [\M:\N]$ (so $\l>\frac{1}{n+2}$).  Moreover, as 
$Jp_i J = p_i$, we may regard the $p_i$'s as elements of ${\rm 
Mat}_m(\mathbb R)$.

The following Lemma \ref{bd} with 
$\epsilon =\frac{1}{2(n+2)}$ then gives
\[
|\mathfrak I(\N,\M)|\leq
(4(n+2)\sqrt{m} + 1)^{m^2}\leq (4(n+2)\sqrt{n} + 1)^{n^2} \ .
\]
\endproof
The following lemma slightly improves \cite[Lemma 2.6]{BD}.
\begin{lemma}\label{bd}
Let $\epsilon>0$ and $\{p_1,p_2,\dots , p_k\}$ be elements in the unit 
ball of ${\rm Mat}_n(\mathbb R)$ such that $||p_i - p_j||\geq\epsilon$ 
if $i\neq j$. Then $k\leq (\frac{2\sqrt n}{\epsilon} + 1)^{n^2}$.
\end{lemma}
\proof
As the uniform and Hilbert-Schmidt norms are related by
$||X||\leq ||X||_{HS}\leq \sqrt{n} ||X||$, the $p_i$'s give
vectors of norm less $\sqrt{n}$ in the Euclidean space $\mathbb R^{n^2}$
(identified with ${\rm Mat}_n(\mathbb R)$ with the Hilbert-Schmidt 
norm) with mutual distance 
larger than $\epsilon$. 
Denoting by $B(r)$ the open ball of radius $r$ in 
$\mathbb R^{n^2}$ we then have
\[
k<\frac{{\rm Vol}(B(\sqrt{n}+\epsilon/2))}{{\rm 
Vol}(B(\epsilon/2))}=(\frac{2\sqrt n}{\epsilon}+1)^{n^2}\ .
\]
\endproof
The following lemmata are variations of known facts (cf. e.g. \cite{T})
and are included for convenience.
\begin{lemma}\label{sp}
Let $x$ be a positive linear operator, $0\leq x\leq 1$, and $p$ 
a selfadjoint projection with $||x-p||\leq\de<1/2$.
Then ${\rm sp}(x)\subset [0,\de]\cup [1-\de,1]$.
\end{lemma}
\proof
With $\ell \in (\de,1 -\de)$ we have $x-\ell 
=(p-\ell)(1+(p-\ell)^{-1}(x-p))$, thus $x-\ell$ is invertible if
$||(p-\ell)^{-1}(x-p)||<1$, which is the case if $\de 
||(p-\ell)^{-1}||<1$. This holds
because $||(p-\ell)^{-1}||=\max\{\ell^{-1},(1-\ell)^{-1}\}<\de^{-1}$. 
\endproof
\begin{lemma}\label{polar}
Let $p$ and $q$ be selfadjoint projections on a Hilbert space $\H$ 
and $t=vh$ be the polar decomposition of $t\equiv pq$. 
If $||p-q||\leq\de<1$, then $v^*v=q$, $vv^*=p$ and $||v-t||\leq\de$.
\end{lemma}
\proof
As $||(p-q)^2||\leq\de^2<1$, the operator $1-(p-q)^2$ is invertible.
Thus $s\equiv pq+(1-p)(1-q)$ is also invertible, indeed
$s^{-1}=(1-(p-q)^2)^{-1} s^*$, and this implies that $v$ is a 
partial isometry from $q$ to $p$. Then we have
\begin{multline*}
||v-t||=||v-vh||= ||v(1-h)||\leq 
||(1-h)\res_{q\H}||
=||(1-\sqrt{qpq})\res_{q\H}||\\
\leq ||(1-qpq)\res_{q\H}||=||q(q-p)q||\leq||q-p||\leq\de\ .
\end{multline*}
\endproof
The bound (\ref{bound}) implies that
\[
|\mathfrak I(\N,\M)|\leq \ell^{\ell} \ ,
\]
where $\ell=[\M:\N]^2$.  The arguments in this section can be 
improved, in particular taking into account that the $p_i$'s are 
projections in Lemma \ref{bd}, leading to a better bound $|\mathfrak 
I(\N,\M)|\leq \ell_1^{\ell_2}$ where however $\ell_2$ is still 
quadratic in the index.  It would be interesting to see if a bound 
$|\mathfrak I(\N,\M)|\leq [\M:\N]^{[\M:\N]}$ holds.  This is the case 
of the example $\N=\M^G$ with $G$ is a finite group where, because of 
the Galois correspondence (see e.g.  \cite{ILP}), $|\mathfrak 
I(\N,\M)|\leq |G|!=[\N:\M]$\, .  We note that 
we have not made use of the specific form of the projections 
associated with intermediate subfactors and the canonical algebra they 
generate \cite{B,BJ}.
\section{Conformal nets and subnets}
We now begin our study of conformal nets. Their subnets will be 
analyzed through the relative superselection structure.
\subsection{Conformal nets}\label{nets}
Let $\I$ denote the family of proper intervals of $S^1$, namely 
connected subsets of $S^1$ of positive measure (length) strictly less than 
$2\pi$. The subnet structure is rather simple for net on the 
4-dimensional Minkowski spacetime (see \cite{CC}), but this does not extend 
to the low-dimensional case.

A net (or precosheaf) $\A$ of von Neumann algebras on $S^1$ is a map 
$$
I\in\I\to\A(I)\subset B(\H)
$$
from $\I$ to von Neumann algebras on fixed a Hilbert space $\H$
that satisfies:
\begin{itemize}
\item[{\bf A.}] {\it Isotony}. If $I_{1}\subset I_{2}$ belong to $\I$, then
\begin{equation*}
 \A(I_{1})\subset\A(I_{2}).
\end{equation*}
\end{itemize}
The net $\A$ is called a (local) {\it conformal} net if in 
addition it satisfies the following properties:
\begin{itemize}
\item[{\bf B.}] {\it Locality}. If $I_{1},I_{2}\in\I$ and $I_1\cap 
I_2=\emptyset$ then 
\begin{equation*}
 [\A(I_{1}),\A(I_{2})]=\{0\},
 \end{equation*}
where brackets denote the commutator\footnote{The locality condition 
will be always assumed in this paper, with the exception of 
Subsection \ref{Fermi}}.
\item[{\bf C.}] {\it Conformal invariance}. 
There exists a strongly 
continuous unitary representation $U$ of $\PSL$ on $\H$ such that
\begin{equation*}
 U(g)\A(I) U(g)^*\ =\ \A(gI),\quad g\in\PSL,\ I\in\I.
\end{equation*}
Here $\PSL$ acts on $S^1$ by M\" obius transformations. We shall 
denote also by $\a_g={\rm Ad}U(g)$ the adjoint action on $B(\H)$.
\item[{\bf D.}] {\it Positivity of the energy}. The generator of the 
one-parameter rotation subgroup of $U$ (conformal Hamiltonian) is positive. 
\item[{\bf E.}] {\it Existence of the vacuum}. There exists a unit 
$U$-invariant vector $\Omega\in\H$ (vacuum vector).
\end{itemize}
We shall say that a conformal net is {\it irreducible} if 
$\vee_{I\in\I}\A(I)=B(\H)$.
Here the lattice symbol $\vee$ denotes the von Neumann algebra 
generated. We recall the following Lemma whose proof can be found in 
\cite{GL2}.
\begin{lemma}\label{irr}
Let $\A$ be a conformal net. The following are equivalent:
\begin{description}
\item[$(i)$] $\A$ is irreducible;
\item[$(ii)$] $\Om$ is cyclic for $\vee_{I\in\I}\A(I)$ and unique $U$-invariant;
\item[$(iii)$] $\Om$ is cyclic for $\vee_{I\in\I}\A(I)$ and the local von Neumann 
algebras $\A(I)$ are factors. In this case they are $III_1$-factors 
(unless $\A(I)=\mathbb C$ identically).
\end{description}
\end{lemma}\noindent
Let $\A$ be an irreducible conformal net. By the Reeh-Schlieder 
theorem \cite{FJ} the vacuum vector $\Om$ is cyclic and separating for each 
$\A(I)$. The {\it Bisognano-Wichmann} property then holds 
\cite{BGL,FG}: the Tomita-Takesaki modular operator $\Delta_I$ and 
conjugation $J_I$ associated with $(\A(I),\Omega)$, $I\in\I$, 
are given by 
\begin{gather}\label{BW} 
U(\Lambda_I (2\pi 
t))=\Delta_{I}^{it},\ t\in\mathbb R,\\ U(r_I)= J_I, 
\end{gather} 
where $\Lambda_I$ is the one-parameter subgroup of $\PSL$ of special 
conformal transformations preserving $I$ and $U(r_I)$ implements a 
geometric action on $\A$ corresponding, the M\" obius reflection on 
$S^1$ mapping $I$ onto $I^\prime$, i.e. fixing the boundary points of 
$I$, see \cite{BGL}.

This immediately implies Haag duality: 
$$
\A(I)'=\A(I'),\quad I\in\I\ ,
$$
where $I'\equiv S^1\setminus I$.
\subsection{Representations}
Let $\A$ be an irreducible local conformal net. A {\it representation} 
$\pi$ of $\A$ is a map
\[
I\in\I\to\pi_I\ ,
\]
where $\pi_I$ is a representation of $\A(I)$ on a fixed Hilbert space 
$\H_{\pi}$ such that
\[ 
\pi_{\tilde I}\res_{\A(I)}=\pi_I, \quad I\subset\tilde I\ ;
\]
we shall always assume that $\pi$ is locally normal, namely
$\pi_I$ is normal 
for all $I\in\I$, which is automatic if $\H_{\pi}$ is separable 
\cite{T}.

We shall say that a representation $\r$ is {\it localized} in a 
interval $I_0$ if $\H_{\r}=\H$ and $\rho_{I'_0}={\rm id}$. Given an 
interval $I_0$ and a representation $\pi$ on a separable Hilbert 
space, there is a representation $\rho$ unitarily equivalent to $\pi$ 
and localized in $I_0$. This is due the type $III$ factor property.

Let $\r$ be a representation of $\A$ localized in a given
interval $I_0$. By Haag duality $\rho$ satisfies the following properties
{\it \begin{description}
\item$(a)$ 
If $I\in\I$ and $I\supset I_0$ then $\rho_I$ is an endomorphism of 
$\A(I)$, and $\rho_{\tilde I}\res_{\A(I)}=\rho_I$ for all $\tilde 
I\in\I$, $\tilde I\supset I$;
\item$(b)$ 
If $I_1\in\I$ and $I_1\cap I_0=\emptyset$, then 
$\rho_{I_1}$ is the identity on $\A(I_1)$;
\item$(c)$ 
If $I,I_1\in\I$ and $I\supset I_0\cup I_1$, there exists a unitary 
$u\in\A(I)$ such that the representation $I\to\rho'_I\equiv 
u\rho_I(\cdot)u^*$ is localized in $I_1$ {\rm (that is to say $\r'_{I_2}$ 
acts identically on $\A(I_2)$ if $I_1\cap I_2 =\emptyset$ for all 
$I_2\in\I$)}.
\end{description}}
We now make a stereographic identification $\mathbb 
R=S^1\setminus\{\infty\}$ and denote by $\I_0\subset\I$ the family of 
bounded intervals of $\mathbb R$, namely of the intervals of $S^1$ whose 
closure do not contain the point $\infty$ of $S^1$.

We denote by $\A_0$ the restriction of $\A$ to $\mathbb R$ (i.e. to 
$\I_0$) and by $\mathfrak A_0$ the associated quasi-local 
C$^*$-algebra $\mathfrak A_0\equiv\overline{\cup_{I\in\I_0}\A(I)}$ 
(norm closure). For a characterization of the so obtained net on 
$\mathbb R$, see \cite{GLW}.

Given $I_0\in\I_0$ a {\it DHR endomorphism} $\rho$ of $\A_0$ localized 
in $I_0$ is a map
\[
\I_0\ni I\to\rho_I
\] 
that associates to each $I\in\I_0$ a representation $\r_I$ of $\A(I)$ 
on $\H$ such that the above conditions $(a), (b), (c)$ hold true with 
$\I$ replaced by $\I_0$.

Clearly a DHR endomorphism determines an endomorphism of $\mathfrak 
A_0$, still denoted by $\r$, such that $\rho_I=\rho\res_{\A(I)}$, 
$I\in\I_0$. The above properties $(a), (b), (c)$ are 
immediately expressed in terms of such endomorphism of $\mathfrak 
A_0$; we shall use the two descriptions interchangeably without 
further specifications.
\begin{proposition}\label{ext}
Let $\rho$ be a DHR endomorphism on $\A_0$ localized in an interval 
$I_0\in\I_0$. There exists a unique representation $\tilde\rho$ of 
$\A$ extending $\rho$ and localized in $I_0$.
\end{proposition}
\proof
Our aim is to define consistently a representation $\rho_I$ of $\A(I)$ 
for every $I\in\I$. To this end, given $I\in\I$, choose $I_1\in\I_0$, 
$I_1\subset I'$ and $L\in\I_0$ let be an interval with $L\supset I_0\cup 
I_1$. Take then a DHR endomorphism of $\A_0$ equivalent to $\rho$ and 
localized in $I_1$, $\rho' ={\rm Ad}u\cdot \rho$ for some unitary 
$u\in\A(L)$. We set
\[
\tilde\rho_I(a) =u^*au, \ a\in\A(I)\ ,
\]
namely $\tilde\rho_I={\rm Ad}u^*\cdot\rho'_I$. Clearly 
$\tilde\rho_I=\rho_I$ if $I\in\I_0$ and a routine checking shows that 
$I\in\I\to\tilde\rho_I$ is indeed a representation of $\A$. 
\endproof 
A representation $\pi$ of $\A$ on a Hilbert space $\H_\pi$ is {\it 
covariant} if there exists a unitary representation $U_\pi$ of the 
universal covering group $\tPSL$ of $\PSL$ on $\H_\pi$ such 
that
\[
{\rm Ad }U_\pi(g)\cdot\pi_I=\pi_{gI}\cdot {\rm Ad } U(g),\quad 
g\in\tPSL,\ I\in\I_0.
\]
Here $U$ has been lifted to $\tPSL$.
$\pi$ is said to have positive energy if the generator of the rotation 
unitary subgroup of $U_\pi$ is positive.

Let $\r$ be a representation of $\A$ localized in $I_0\in\I$.  By a 
{\it local cocycle} (w.r.t.  to $\r$) we shall mean the assignement of 
an interval $I\supset I_0$, a symmetric neighborhood $\U$ of the 
identity of $\tPSL$ such that $I_0\cup gI_0\subset I, \forall g\in\U$ 
and a strongly continuous unitary valued map 
$z: g\in\U\to z_\r(g)\in\A(I)$ such that
\begin{gather}\label{cov1} z_\r(g)\in\A(I)\\ 
z_\r(gh)= z_\r(g)\a_g(z_\r(h)),\label{cov2}\\ 
{\rm Ad } z_\r(g)^*\cdot\r_{\tilde 
I}(a)= \a_g\cdot\r_{g^{-1}\tilde I}\cdot\a_{g^{-1}}(a),
\ a\in\A(\tilde I),\label{cov3} 
\end{gather}
for some open interval $\tilde I$ with $\bar I\subset\tilde I$
and all $g,h\in\U$ such that $I\cup gI\subset\tilde I$.   
We shall then say that $z$ is localized in $I$.

If this holds, then eq. (\ref{cov3}) is valid for all $L\in\I$:
\begin{equation}\label{global}
{\rm Ad } 
z_\r(g)^*\cdot\r_L(a)=\a_g\cdot\r_{g^{-1}L}\cdot\a_{g^{-1}}(a),\quad
a\in\A(L).
\end{equation}
Indeed, if $L\supset I$ then the above equation holds
by additivity \cite{FJ}. Thus it holds for sub-intervals $L_0\subset L$. 
Again by additivity, the equation is then satisfied for all $L\in\I$.

If $\r$ is a covariant representation of $\A$ localized in $I_0$ then 
for any given interval $I\supset I_0$ there exists a local 
cocycle w.r.t. $\r$ localized in $I$.
Indeed if $\U$ is a symmetric
neighborhood of the identity of $\tPSL$ such that $I_0\cup gI_0\subset I, 
\forall g\in\U$, then by Haag duality the unitaries
\[
z_\r(g)\equiv U_\r (g)U(g)^*	
\]
belong to $\A(I)$ for all $g\in \U$ and 
clearly verify the local cocycle property (\ref{cov2},\ref{cov3}). 

Notice now that, taking $I_0,I\in\I_0$, a local cocycle is expressed 
in terms of the DHR endomorphism of $\A_0$ associated with $\r$. The 
converse construction is made in the following.

\begin{proposition}
Let $\rho$ be a DHR endomorphism of $\A_0$ localized in the interval 
$I_0\in\I_0$ and $\tilde\rho$ the representation of $\A$ extending 
$\rho$ given by Proposition \ref{ext}. Then $\tilde\rho$ is covariant 
iff there exists a local cocycle $z_{\r}$ w.r.t. $\r$ (i.e. 
properties (\ref{cov1},\ref{cov2},\ref{cov3}) hold with $\I$ replaced 
by $\I_0$).
\end{proposition}
\proof
We need only to show that $\tilde\rho$ is covariant if there exists a 
local cocycle $z_{\r}$. By the above arguments eq. (\ref{global}) holds.
Now set $U_{\tilde\r}(g)=z_\r(g)U(g)$ for 
$g$ in a suitable neighborhood of the identity of $\tPSL$. Then 
$U_{\tilde\r}$ is a local representation of $\tPSL$, hence it extends 
to a unitary representation of $\tPSL$ because $\tPSL$ is simply 
connected. The local covariance then gives
\[
U_\r(g)\r_L(a)U^*_\r(g)=\r_{gL}(\a_g(a)),\ a\in\A(L),
\]
for any $L\in\I$. The covariance then follows by the group property 
of $U_\r$, see also \cite{GL1}.
\endproof
Before concluding this section we recall that, if $\r$ is a localizable 
representation of $\A$, the (statistical) dimension of $\r$ is 
$d(\r)\equiv[\r_{I'}(\A(I'))' : \r_I(\A(I))]^{\frac{1}{2}}$, 
independently on $I\in\I$, and this clearly coincides with 
$[\A(I) : \r_I(\A(I))]^{\frac{1}{2}}$ if $\r$ is localized in $I$
\cite{L2}.

If $\r$ is M\"{o}bius covariant and $d(\r)<\infty$, then $\r$ has
positive energy \cite{BCL}.
\subsection{Subnets}
Let $\A$ be a local irreducible conformal net of von Neumann algebras 
on $S^1$ as above and $U$ the associated unitary positive energy 
representation of $\PSL$ on the vacuum Hilbert space $\H$.

By a \emph{conformal subnet} we shall mean a map 
\[
I\in\I\to\B(I)\subset\A(I)
\] 
that associates to each interval $I\in\I$ a von Neumann subalgebra $\B(I)$ 
of $\A(I)$, which is isotone
\[
\B(I_1)\subset\B(I_2),\quad I_1\subset I_2 \ ,
\]
and M\" obius covariant w.r.t. the representation $U$, namely
\[
U(g)\B(I)U(g)^{-1}= \B(gI)
\]
for all $g\in \PSL$ and $I\in\I$.

Let $\H_{\B}$ be the closure of $(\vee_{I\in\I}\B(I))\Om$ and $E$ the 
orthogonal projection of $\H$ onto $\H_{\B}$. By the Reeh-Schlieder 
theorem $\overline{\B(I)\Om}=\H_{\B}$ for each fixed $I\in\I$. 

Clearly $\H_\B$ is $U$-invariant and $\Om$ is unique 
$U\res_{\H_{\B}}$-invariant, thus by Lemma \ref{irr} the restriction of $\B$ to $\H_{\B}$ 
is an irreducible local conformal net on $\H_\B$ where $U\res_{\H_\B}$ 
is the associated unitary representation of $\PSL$.

As $\Omega$ is separating for $\A(I)$, $\Omega$ is also separating for 
$\B(I)$, $I\in\I$. Thus the restriction map $b\in\B(I)\to 
b\res_{\H_\B}$ is is one-to-one, so we will often identify $\B$ with 
its restriction to $\H_\B$; should we need to specify, we shall talk 
on the net $\B$ on $\H$ or on $\H_\B$. Note that each $\B(I)$ is a 
factor.
\begin{lemma} 
For each $I\in\I$ there is a vacuum preserving conditional expectation 
$\e_I:\A(I)\to\B(I)$ such that $\e_{\tilde I}\res_{\A(I)}=\e_I$ if 
$I\subset\tilde I$. Thus $\B$ is a standard net of subfactors in the 
sense of \cite{LR}.
\end{lemma}
\proof
By the Bisognano-Wichmann property $\B(I)$ is globally invariant 
under the modular group of $(\A(I),\Om)$, hence by Takesaki's theorem 
there exists a conditional expectation $\e_I:\A(I)\to\B(I)$ given by
\[
\e_I(a)E = EaE,\: a\in\A(I)\ . 
\]
As $E$ is independent of $I$, we have that $\e_{\tilde I}\res_{\A(I)}=\e_I$ if 
$I\subset\tilde I$.
\endproof
By M\" obius covariance the index $[\A(I):\B(I)]$ is independent of 
the interval $I\in\I$ and will be denoted by $[\A:\B]$. The 
following lemma is contained in \cite{BE} with the strong additivity 
assumption and in \cite{DLR} in the conformal case.
\begin{lemma} If $[\A:\B]<\infty$ then $\B(I)'\cap\A(I)=\mathbb C$, 
$I\in\I$.
\end{lemma}
\proof
By the Bisognano-Wichmann property and the uniqueness of the vacuum 
the modular group of $\A(I)$ w.r.t. $\Om$ acts ergodically on 
$\B(I)'\cap\A(I)$, hence $\B(I)'\cap\A(I)=\mathbb C$ because 
$\B(I)'\cap\A(I)$ is finite-dimensional.
\endproof
We shall make a variation of the analysis made in \cite{LR}, which is 
needed because our nets are not directed.
\begin{lemma}\label{isometry}
Let $[\A:\B]<\infty$ and $I_0\subset I$. There exists
a canonical endomorphism $\g_{ I}:\A(I)\to\B(I)$ 
with associated Q-system $(\g_{I},T,S)$ such that $T\in \A(I_0)$ 
and $S\in\B(I_0)$ and $\g_{I}\res_{\A(I_0)}$ 
is a canonical endomorphism of $\A(I_0)$ into $\B(I_0)$.
\end{lemma}
\proof Let $\C(I_0)=\langle \A(I_0),E\rangle$ and $\C(I)=
\langle \A(I),E\rangle$ be the Jones extensions and 
$\e'_{I_0}:\C(I_0)\to\A(I_0)$, $\e'_{I}:\C(I)\to\A(I)$ 
the dual expectations. Since every $X\in\C(I_0)$ can be written 
as $X=\sum_i x_i Ey_i$ with $x_i,y_i\in\A(I_0)$ \cite{J,PP}, we then have 
\[
\e'_{I}(X)=\sum_i x_i\e'_{I}(E)y_i
=\l\sum_i x_i y_i=\sum_i x_i\e'_{I_0}(E)y_i = \e'_{I_0}(X),\ X\in\A(I),
\]
where $\l\equiv [\A:\B]^{-1}$, namely $\e'_{I}\res_{\A(I_0)}=\e'_{I_0}$.

Let $V\in \C(I_0)$ be an isometry $VV^*=E$. Then a canonical endomorphism 
$\g_{I}:\A(I)\to\B(I)$ is given by 
\begin{equation}
\g_{I}(a)E=VaV^*\ .	
	\label{vav}
\end{equation}
Now 
\[
T=\l^{-1}\e'_{I}(V),\quad S=\l^{-1}\e_{I}(T),
\]
are the isometries in the Q-system for $\g_{I}$ and $T\in 
\A(I_0)$ and $S\in\B(I_0)$ by the compatibility of the expectations.
\endproof
Notice the formula
\begin{equation}\label{e-can}
\g_I(a)=\l^{-1}\e'_I(VaV^*), \quad a\in\A(I) ,
\end{equation}
which is obtained applying $\e'_I$ to both members of eq. (\ref{vav}).
\begin{proposition}\label{rep-can}
Let $\A$ be a local irreducible conformal net on $S^1$ and 
$\B\subset\A$ a conformal subnet. Given an interval $I_0\in\I$ the 
dual canonical endomorphism $\theta_{I_0}\equiv\g_{I_0}\res_{\B(I_0)}$ 
extends to a representation $\theta$ of $\B$ localized in $I_0$.
\end{proposition}
\proof
First we assume $[\A:\B]<\infty$.  By Proposition \ref{ext} it is 
sufficient to show that $\theta_{I_0}$ extends to a DHR endomorphism 
of $\B$ localized in $I_0$ (properties $(a),(b),(c)$ with $\I$ 
replaced by $\I_0$).  This is soon verified by applying Lemma 
\ref{isometry}.  Indeed property $(a)$ is an immediate consequence of 
this lemma.  Concerning property $(b)$ notice the formula \cite{L3}
\[
\theta_I(b)=\l^{-1}\e_I(TbT^*), \ b\in\B(I),
\]
with $T$ as in Lemma \ref{isometry}, that follows similarly to the 
formula (\ref{e-can}). Then $(b)$ follows because $T\in\A(I_0)$.

Finally property $(c)$ is immediate by the uniqueness up to inners of the 
canonical endomorphism \cite{L}.

The general case can be obtained along the same lines making use of 
\cite[Theorem 3.2]{LR} instead of Lemma \ref{isometry}.
\endproof
\begin{proposition}\label{thetacovariance}
The representation $\theta$ of $\B$ on $\H_\B$ in Prop. \ref{rep-can} 
is unitarily equivalent to the identity representation of $\B$ on 
$\H$. In particular $\theta$ is covariant with positive energy.
\end{proposition}
\proof
Indeed the isometry $V\in\C(I)$ with $VV^*=E$ (as in Lemma 
\ref{isometry}) satisfies the equation 
\begin{equation}\label{v}
\theta_{\tilde I}(b)E=VbV^*,\ b\in \B(\tilde I)\ ,
\end{equation}
if $\tilde I$ is an interval containing $I$.

To show that $V$ implements the desired unitary equivalence we need to 
further show that the above equation (\ref{v}) holds true with $\tilde 
I$ replaced with an interval $I_1$ not containing $I$.  This is 
certainly true if $I_1\cap I=\emptyset$, because in this case $\theta$ 
acts trivially on $\B(I_1)$ and $V\in\C(I)$ commutes with $\B(I_1)$ 
because $\C(I)=\langle\A(I),E\rangle$ and $\B(I_1)$ commute.

So we may assume that $I_1\supset I'$, extending $I_1$ if necessary.  
Choose then an interval $I_0\subset I$ with $I_0\cap I_1=\emptyset$.  
By Lemma \ref{isometry} we can find a canonical endomorphism 
$\A(I_0)\to\B(I_0)$ with dual canonical endomorphism extending to a 
representation $\theta'$ of $\A$ localized in $I_0$ with a unitary 
$u\in\B(I)$ such that $\theta'={\rm Ad}u\cdot\theta$.  Then the 
isometry $V'$ associated with $\theta'$ belongs to $\C(I_0)$ and is 
given $V'=uV$.  Therefore
\[
\theta(b)E=u\theta'(b)u^*=uV'b{V'}^*u^* =VbV^*, \ b\in \B(I_1)
\]
as desired.
\endproof
The following corollary is a consequence of the equivalence between 
local and global intertwiners for a finite-index covariant 
representation \cite{GL2}. 
\begin{corollary}\label{inter}
Assume $[\A:\B]<\infty$ and let $\theta$ be the representation of $\B$ 
on $\H_{\B}$ in Proposition \ref{rep-can}.  Then $\theta$ has a finite 
direct sum decomposition 
\begin{equation}\label{dec-rep} 
\theta=\bigoplus_{i=0}^N N_i \r_i\ ,
\end{equation}
where the $\r_i$'s are positive-energy covariant irreducible 
representation of $\B$ on $\H_\B$ localized in $I_0$. Thus 
$\theta_I=\oplus_{i=0}^N N_i \r_{iI}$ is a decomposition of the 
canonical endomorphism for any interval $I\supset I_0$.
\end{corollary}
\proof
As $d(\theta)=[\A:\B]<\infty$ we may decompose $\theta$ into irreducible 
representations $\r_i$ as above localized in $I_0$; moreover the 
$\r_i$ are covariant because $\theta$ is covariant, see \cite{GL2}.
Then by \cite[Theorem 2.3]{GL2} each $\r_{iI}$ is an irreducible 
endomorphism of $\B(I)$ if $I\supset I_0$.
\endproof
Thus, if $[\A:\B]<\infty$ then the identity representation of $\B$ on 
$\H$ has finite statistical dimension.  The converse is also true: as 
we have the inclusions $\B(I)\subset\A(I)\subset\B(I')'$, if the 
identity representation of $\B$ on $\H$ has finite statistical 
dimension, namely $[\B(I')':\B(I)]<\infty$, then $[\A(I):\B(I)]<\infty$.
\begin{corollary}\label{zg}
Assume $[\A:\B]<\infty$.  With the above notations, let 
$K_i\subset\A(I_0)$ be the Hilbert spaces of isometries corresponding 
to $\r_{iI_0}$ as in eq.  (\ref{Ki}) for the inclusion 
$\B(I_0)\subset\A(I_0)$.  If $I$ is an interval and $I\supset I_0$, 
then $\A(I)$ is generated by $\B(I)$ and $K_i$ as in 
(\ref{expansion}).

If $R_i\in K_i$ then $\a_g(R_i)=z_{\r_i}(g)^* R_i$, for all 
$g\in\PSL{\tilde{\phantom x}}$ such that $gI_0\subset I$, where 
$z_{\r_i}(g)\in\B(I)$ are unitaries in the local cocycle associated 
with $\r_i$ (\ref{cov1}).
\end{corollary}
\proof
By Corollary \ref{inter} we have 
$(\r_{iI},\theta_I)=(\r_{iI_0},\theta_{I_0})\subset\B(I_0)$, therefore 
by formula (\ref{frob}) and Lemma \ref{isometry} $K_i\subset\A(I_0)$ 
is also the Hilbert space associated with $\B(I)\subset\A(I)$ as in 
(\ref{expansion}). By Lemma \ref{lemmaexpansion} $\A(I)$ is then 
generated by $\B(I)$ and the $K_i$'s.

Let $g\in\tPSL$ such that $gI_0\subset I$. Then 
$\a_g\cdot\r_i\cdot\a_{g^{-1}}$ is localized in $I$ and 
$\a_g\cdot\r_i\cdot\a_{g^{-1}}={\rm Ad}z_{\r_i}(g)^*\cdot\r_i$ by 
formula (\ref{cov3}) and Cor. \ref{inter}, where $z_{\r_i}(g)\in\B(I)$. 
Therefore $\a_g(K_i)=z_{\r_i}(g)^* K_i$, namely 
$\a_g(R_i)=z_{\r_i}(g)^* D(g)^* R_i$ for all $R_i\in K_i$, where 
$D(g)\in B(K_i)$. It is immediate to check that $D$ locally satisfies 
the cocycle property with respect to $\a^{\r_i}_g\equiv {\rm 
Ad}z_{\r_i}(g)\cdot\a_g$, namely $D(gh)=D(g)\a^{\r_i}_g(D(h))$ for 
$g,h$ in a suitable neighborhood of $\tPSL$. But 
$\a^{\r_i}\res_{B(K_i)}$ is a finite-dimensional representation of 
$\tPSL$, thus it must be trivial because $\tPSL$ has no non-trivial 
unitary finite-dimensional representations. Thus $D$ is a local 
finite-dimensional unitary representation of $\tPSL$ on $K_i$, so $D$ 
is again trivial.

Therefore $\a_g(R_i)=z_{\r_i}(g)^* R_i$, for all 
$g\in\PSL{\tilde{\phantom x}}$ such that $gI_0\subset I$ as desired.
\endproof
If $\A$ is a conformal net on a Hilbert space $\H$ and $\B$ is a 
conformal subnet, we shall set $\B'\equiv(\vee_{L\in\I}\B(L))'$, where 
the commutants are taken on $\H$.  For completeness we mention the 
following partial extension of Cor.  \ref{inter} to the infinite index 
case, although it is not used in this paper. 
\begin{proposition}\label{univac}
Let $\A$ be a local irreducible conformal net on $S^1$ and 
$\B\subset\A$ a strongly additive 
conformal subnet. The following are equivalent: 
\begin{description}
\item[$(i)$] The identity representation of $\B$ on $\H$ contains the 
vacuum representation of $\B$ with multiplicity one.
\item[$(ii)$] The identity representation of $\B$ on $\H$ contains the 
vacuum representation of $\B$ with finite multiplicity.  
\item[$(iii)$] $\B^c=\mathbb C$, where
$\B^c(I)\equiv\B'\cap\A(I)$.
\end{description}
\end{proposition}
\proof 
$(iii)\implies (i)$.  By Proposition \ref{thetacovariance} we have to 
show that the intertwiner space between the representation $\theta$ on 
$\H_\B$ and the identity representation of $\B$ on $\H_\B$ is 
one-dimensional.  If $\theta$ is localized in $I$ then, by Haag 
duality, any such intertwiner belongs $(\theta_I,\iota)$ namely it 
belongs to $\B(I)$ and intertwines $\theta_I$ and the identity 
automorphism of $\B(I)$.  But $(\theta_I,\iota)$ is one-dimensional 
because $\B(I)\subset\A(I)$ is an irreducible inclusion of infinite 
factors with a normal conditional expectation \cite{L2,FI}.

$(i)\implies (ii)$ is obvious, we show $(ii)\implies (iii)$.  Denote 
by $\pi$ the subrepresentation of the identity representation of $\B$ 
on $\H$ corresponding to the vacuum representation and $\K\subset\H$ 
the corresponding subspace.  Then we have a decomposition 
$\K=\H_0\otimes\H'_0$ and $\pi=\pi_0\otimes{\rm id}$ where $\pi_0$ 
is irreducible and ${\rm dim}\H'_0<\infty$ The representation $U$ of 
$\PSL$ decomposes as $U=U_0\otimes U'_0$ and, since $\PSL$ has no 
non-trivial unitary representation, $U'_0$ is the identity on $\H'_0$.  
As $\Om$ is unique $U$-invariant, we then have ${\rm dim}\H'_0=1$, 
namely $(i)$ holds.

Thus $\K=\H_\B$ and the projection $E$ onto $\H_\B$ belongs to the 
center of $\B'$.  In particular if $b\in\B^c(I)$ then $b$ commutes 
with $E$.  But $E$ implements the expectation $\e_I$ and $b\in\A(I)$, 
so $b$ belongs to the center of $\B(I)$, thus $b$ is a scalar.
\endproof
\subsection{Intermediate subnets}
Let $\A$ be a local irreducible conformal net and $\B$ a conformal 
subnet with finite index. We now show that there exists a one-to-one 
correspondence between $\mathfrak I(\B(I_0),\A(I_0))$ for a fixed 
interval $I_0$ and the set of intermediate conformal nets between 
$\B$ and $\A$.
\begin{theorem}\label{int-net}
Let $I_0$ be a fixed interval of $S^1$ and $\R$ be an intermediate 
subfactor between $\B(I_0)$ and $\A(I_0)$.
There exists a unique conformal subnet $\L$ on $\A$ with 
$\B(I)\subset\L(I)\subset\A(I)$ and $\L(I_0)=\R$.
\end{theorem}
\proof
Let $\Lambda_I$ denote as before the one-parameter group of 
special conformal transformations preserving $I$. As is easily seen 
$\Lambda_I(\mathbb R)$ is exactly the subgroup of $\PSL$ of those $g$ 
with $gI=I$. Then $t\to\b_t ={\rm Ad } 
U(\Lambda_{I_0}(t))\res_{\A(I_0)}$ is a one-parameter automorphism 
group of $\A(I_0)$ leaving $\B(I_0)$ globally invariant, hence 
$\b_t(\R)=\R$ by Corollary \ref{invar}.

Now, for $I\in\I$ we set 
\[
\L(I)\equiv \a_g(\R)\ ,
\] 
where $g\in \PSL$ is a M\" obius transformation such that $I_0 = gI$
and $\a_g\equiv{\rm Ad } U(g)$.
$\L(I)$ is indeed 
well-defined because if $h\in \PSL$ is any other element with $I_0 = 
hI$, then $h^{-1}g I_0=I_0$, thus $h^{-1} g=\Lambda_{I_0}(t)$ for some 
$t\in\mathbb R$ and $\a_{h^{-1}g}(\R)=\b_t(\R)=\R$, thus 
$\a_g(\R)=\a_h(\R)$.

As $\L(I)\subset\A(I)$ and $\A$ is local, $\L(I_1)$ and $\L(I_2)$ 
clearly commute if $I_1\cap I_2=\emptyset$. To 
show that $I\in\I\to\L(I)$ is a conformal net we need to check the 
isotony property, namely that $\L(I_1)\subset\L(I)$ if $I_1\subset 
I$ are intervals. By conformal 
invariance we may assume that $I_1=I_0$ and that $I=g_0 I_0$ for some
$g_0\in \PSL$ and then we need to show that 
$\a_{g_0}(\R)\supset\R$. 

Now by Corollary \ref{zg} $\A(I)$ is generated by $\B(I)$ and Hilbert 
spaces of isometries 
$K_i\in\A(I_0)$ corresponding to the expansion (\ref{expansion}). 
Moreover 
\[
\a_{g_0}(R_i) = z_{\r_i}(g_0)^* R_i ,\quad \forall R_i\in K_i\ ,
\]
where $z_{\r_i}(g_0)\in\B(I)$ by Corollary \ref{zg}. With $K'_i\equiv 
K_i\cap\R$, by Lemma \ref{zg} $\R$ is generated by $\B(I_0)$ and the 
$K'_i$'s, thus
\begin{multline*}
\a_{g_0}(\R)=\{\a_{g_0}(\B(I_0)),\a_{g_0}(K'_i)\}''
=\{\B(I),z_{\r_i}(g_0)^* K'_i\}''\\
=\{\B(I),K'_i\}''\supset\{\B(I_0),K'_i\}''=\R
\end{multline*}
as desired.
\endproof
\subsection{Complete rationality is hereditary}
\label{Hered}
To simplify notations, given two different points ${\bf a},{\bf b}$ of 
$S^1$, we shall write $[{\bf a},{\bf b}]$ for the closure of the set 
of all $z\in S^1$ that follow ${\bf a}$ and precede ${\bf b}$ in the 
counterclockwise order, and by $({\bf a},{\bf b})$ the interior of 
$[{\bf a},{\bf b}]$.  Two intervals $I_1$, $I_2$ will be called {\it 
adjacent} if there are three different points ${\bf a},{\bf b},{\bf 
c}\in S^1$ such that $\bar I_1 = [{\bf a},{\bf b}]$, $\bar 
I_2 = [{\bf b},{\bf c}]$ and $\bar I_1\cup \bar I_2$ belongs to $\I$.

Let $\A$ be a conformal net on $S^1$. Recall that $\A$ is {\it 
strongly additive} if
\[
\A(I_1)\vee\A(I_2)=\A(I)\ ,
\]
if $I_1$, $I_2$ are adjacent intervals and $I=\bar I_1\cup\bar 
I_2\in\I$.

$\A$ is {\it split} if $\A(I_0)\subset\A(I)$ is a split inclusion of 
von Neumann algebras, namely $\A(I_0)\vee\A(I)'$ is naturally 
isomorphic to $\A(I_0)\otimes\A(I)'$, if $I_0, I\in\I$ and $\bar I_0$ 
contained in the interior of $I$.

If $\A$ is split and $I_1,I_2$ are intervals with disjoint closures, 
then $\A(I_1)\vee\A(I_2)\simeq\A(I_1)\otimes\A(I_2)$ is a factor and 
we shall denote by $\mu_\A$ the index of the 2-interval inclusion 
$\A(I_1)\vee\A(I_2)\subset (\A(I_3)\vee\A(I_4))'$ where
$I_3,I_4$ are the two connected components of $S^1\setminus (I_1\cup 
I_2)$.

We shall say that $\A$ is {\it completely rational} if $\A$ is split, strongly additive 
and the index $\mu_\A<\infty$ where the $I_i$'s are intervals as above.
\begin{lemma}\label{usplit} 
Let $\A$ be an irreducible local conformal net on $S^1$ and $\B\subset 
\A$ a finite-index conformal subnet.  Then $\A$ is 
split and $\mu_\A<\infty$ iff $\B$ is split and $\mu_\B<\infty$.

In this case the relation
$\mu_\B =[\A:\B]^2\mu_\A$ holds.
\end{lemma}
\proof
First notice that, if $\A$ and $\B$ are split, the proof of 
\cite[Prop. 24]{KLM} shows that $[\A:\B]^2\mu_\B =[\A:\B]^4\mu_\A$, 
thus $\mu_\B =[\A:\B]^2\mu_\A$ because $[\A:\B]<\infty$.
So the lemma is proved if we show the following implications:
\begin{align}\label{a}
\A\ \text{split}\ \&\ \mu_\A<\infty &\implies \B\  \text{split}\ ,\\
\B\ \text{split}\ \&\ \mu_\B<\infty &\implies \A\ \text{split}\ 
.\label{b}
\end{align}
Now the implication $\A$ split $\implies$ $\B$ split is rather 
immediate so, by the above comment, the first implication (\ref{a}) holds.

The second implication (\ref{b}) is proved in \cite[Prop. 25]{KLM} 
in a specific case, but the 
argument given there works in general. For the convenience of the 
reader we make this explicit.

Let $I_1, I_2$ be intervals with disjoint closures and $I_3,I_4$ the bounded connected
components of $I'_1\cap I'_2$. The conditional
expectation  $\e_I:\A(I)\to\B(I)$ associated with the interval $I$,
where $I$ is the interior of
$\bar I_1\cup \bar I_2\cup \bar I_3$, maps $\A(I_1)\vee\A(I_2)$ onto
$\B(I_3)'\cap\B(I)=(\B(I_3)\vee\B(I_4))'$, thus
\begin{equation*}
\e\equiv\e_0\cdot\e_I |_{\A(I_1)\vee\A(I_2)}
\end{equation*}
is a normal faithful expectation of $\A(I_1)\vee\A(I_2)$ onto 
$\B(I_1)\vee\B(I_2)$,
where $\e_0$ is a normal faithful expectation of $(\B(I_3)\vee\B(I_4))'$ onto 
$\B(I_1)\vee\B(I_2)$, that exists because $\mu_{\B}<\infty$.

To get the split property of $\A$,
it will suffice 
to show that the above expectation $\e$ satisfies
\[
\e(a_1 a_2) = \e(a_1)\e(a_2) \ ,\quad a_i\in\A(I_i) \ ,
\]
and $\e(\A(I_i))\subset \B(I_i)$, as we may then compose a normal
product state
$\phi_1\otimes\phi_2$ of $\B(I_1)\vee\B(I_2)\simeq
\B(I_1)\otimes\B(I_2)$ with $\e$ to get a normal product state of
$\A(I_1)\vee\A(I_2)$.

Let $R_{i,k}^{(\ell)}\in\A(I_{\ell})$, $\ell=1,2$, be elements satisfying
the relations (\ref{expansion}) for the inclusion
$\B(I_i)\subset\A(I_i)$, so that $\A(I_{\ell})$ is generated by 
$\B(I_{\ell})$
and the $R_{i,k}^{(\ell)}$'s. With $a^{(\ell)}\in\A(I_{\ell})$ we then 
have an expansion
\[
a^{(\ell)}=\sum_{i,k} b_{i,k}^{(\ell)} R_{i,k}^{(\ell)} \ ,\quad b_{i,k}^{(\ell)} 
\in\A(I_{\ell})\ ,
\]
hence
\[
a^{(1)}a^{(2)}=\sum_{i,h,j,k} b_{i,h}^{(1)} b_{j,k}^{(2)} R_{i,h}^{(1)} R_{j,k}^{(2)}\ ,
\]
so we have to show that $\e(R_{i,h}^{(1)} R_{j,k}^{(2)})=0$ unless
$i=j=0$. Now $R_{i,h}^{(1)}=u_{i,h} R_{i,h}^{(2)}$ for some unitary
$u_{i,h}\in (\B(I_3)\vee\B(I_4))'\subset\B(I)$ 
and we have
\begin{multline*}
\e (R_{i,h}^{(1)} R_{j,k}^{(2)})=\e (u_{i,h} R_{i,h}^{(2)} R_{j,k}^{(2)})
=\e_0 (u_{i,h} \e_I (R_{i,h}^{(2)} R_{j,k}^{(2)}))\\
=\e_0 (u_{i,h} \e_{I_2} (R_{i,h}^{(2)} R_{j,k}^{(2)}))
= \e_0 (u_{i,h})\e_{I_2} (R_{i,h}^{(2)} R_{j,k}^{(2)}) \ .
\end{multline*}
As $\e_0 (u_{i,h})\in\B(I_1)\vee\B(I_2)$ is        
an intertwiner between irreducible  endomorphisms localized
in $I_1$ and $I_2$, we have $\e_0 (u_{i,h})=0$, thus $\e (R_{i,h}^{(1)} 
R_{j,k}^{(2)})=0$, for all $i\neq 0$.
If $i=0$ and $j\neq 0$, then again $\e (R_{i,h}^{(1)} R_{j,k}^{(2)})=
\e (R_{j,k}^{(2)})=0$ because $\e_I (R_{j,k}^{(2)})=\e_{I_2}(R_{j,k}^{(2)})=0$.
\endproof
\begin{lemma}\label{uss} 
Let $\A$ be an irreducible local conformal net on $S^1$ and $\B\subset 
\A$ a finite-index conformal subnet. Then $\A$ is strongly additive 
if $\B$ is strongly additive.
\end{lemma}
\proof 
Let $I_1,I_2$ be adjacent intervals with 
$I\equiv{\bar I_1\cup \bar I_2}$ and let $T\in\A(I_1)$
be the isometry in the Q-system for $\g_{I_1}$ as in Lemma 
\ref{isometry}. Then, by applying Lemma 
\ref{isometry}, $T$ is also the isometry in the Q-system associated 
with $\g_{I}$. 
In particular $\A(I_1)=\B(I_1)T$ and $\A(I)=\B(I)T$, thus
\[
\A(I_1)\vee\A(I_2)=\B(I_1)T\vee\A(I_2)\supset 
\{\B(I_1)\vee\B(I_2),T\}''=\{\B(I),T\}''=\A(I)\ .
\]
This concludes the proof.
\endproof 
\begin{theorem}\label{str.add.} 
Let $\A$ be an irreducible local conformal net on $S^1$ and $\B\subset 
\A$ a conformal subnet with $[\A:\B]<\infty$.
Then $\A$ is completely rational iff $\B$ is completely rational.
\end{theorem}
\noindent
We postpone the proof of this theorem to Subsection \ref{proof}.
\subsubsection{Complete rationality of coset models}
\label{list}
We begin with the following simple lemma.
\begin{lemma}\label{co}
Let $\A_1$ and $\A_2$ be irreducible local conformal nets on $S^1$. 
Then $\A_1\otimes\A_2$ is split (resp. strongly additive; completely 
rational) iff both $\A_1$ and $\A_2$ are split (resp. strongly 
additive; completely rational).
\end{lemma}
\proof
All this can be checked directly, see \cite{DL,KLM}.
\endproof 
Let $\A$ be an irreducible local conformal net on $S^1$ and $\B\subset 
\A$ a conformal subnet.  Then 
\[
\B^c:I\in\I\to\B^c(I)\equiv\B'\cap\A(I)
\]
is clearly a conformal subnet of $\A$ and is called the {\it coset} 
net associated with $\B\subset\A$, cf. \cite{X1}.  Also $\B\vee\B^c 
:I\in\I\to\B(I)\vee\B^c(I)$ is then a conformal subnet of $\A$.

Following F. Xu \cite{X1}, we call $\B$ {\it cofinite} in $\A$ if 
$[\A:\B\vee\B^c]<\infty$.  Indeed in \cite{X1} $\B$ is strongly 
additive so $\B^c(I)=\B(I)'\cap\A(I)$ is the relative commutant of 
$\B(I)$ in $\A(I)$.

By Takesaki theorem $\B(I)\vee\B^c(I)$ is naturally isomorphic to the 
von Neumann tensor product $\B(I)\otimes\B^c(I)$.  
\begin{corollary}\label{cofinite} 
Let $\A$ be an irreducible local conformal net on $S^1$ and $\B\subset 
\A$ a cofinite conformal subnet.  With the above notations, $\A$ is 
completely rational iff both $\B$ 
and $\B^c$ are completely rational.

If $\A$ is `split \& strongly additive', so are $\B$ and $\B^c$.
\end{corollary}
\proof
It is enough to apply Theorem \ref{str.add.} and 
Lemma \ref{co} and Proposition \ref{ss} below.
\endproof
To give a first application, suppose now that the net $\A$ is 
diffeomorphism invariant.  Then one can consider the conformal subnet 
$\A_{{\rm Vir}(c)}\subset\A$, which is associated with the vacuum 
representation of the Virasoro algebra with central charge $c>0$ see 
e.g.  \cite{BS}.
\begin{corollary}
Let $\A$ be split, strongly additive and diffeomorphism invariant. If
$\A_{{\rm Vir}(c)}$ is cofinite in $\A$, then $c\leq1$.
\end{corollary}
\proof If $\A_{{\rm Vir}(c)}$ is cofinite in $\A$ then by Cor.  
\ref{cofinite} also $\A_{{\rm Vir}(c)}$ is strongly additive, which is 
not the case if $c>1$ \cite{BS}.  
\endproof 
We now turn our attention to coset models.  Let $G$ be simply 
connected semisimple compact Lie group of type $A$, i.e.  
$G=SU(N_1)\times SU(N_2)\times\cdots \times SU(N_n)$.  If $H\subset G$ 
is a closed subgroup, there is a corresponding inclusion of loop 
groups $LH\subset LG$.  Then the vacuum representation of $LG$ at 
level $k$ (see \cite{PS}) gives an inclusion of conformal nets denoted 
by $H\subset G_k$ (where $H$ may also have a suffix denoting the 
appearing level).  Thanks to results of Xu \cite[Corollary 3.1]{X1}, 
see also the correction in \cite{X3}, the inclusions of conformal nets 
associated with
\begin{description}
\item[$(i)$] $G_{k_1 + k_2 + \cdots +k_m}\subset 
G_{k_1}\times G_{k_2}\times\cdots G_{k_m}$ where the inclusion is 
diagonal, $k_i\in\mathbb N$, $i= 1,\dots,m$ and $G=SU(n)$,

\item[$(ii)$] $H_{\ell k}\subset G_{\ell}$, if $H_{k}\subset G_{1}$ is a 
conformal inclusion, where $k$ is the Dynkin index, $\ell\in\mathbb 
N$, $H$ is simple and of type $A$ and $G$ is simple,

\item[$(iii)$] $H\subset G_m$, where $H$ is the Cartan subgroup of $G$,
\end{description}
are all cofinite. So we have the following corollary.
\begin{corollary}
The coset subnets corresponding to the inclusions of nets in $(i)$, 
$(ii)$, $(iii)$ of the above Xu's list are completely rational.
\end{corollary}
\proof
As the conformal net $SU(N)_k$ is completely rational \cite{X5} (a 
correct proof of the strong additivity is contained in \cite{To}), and 
the subnets in Xu's list are cofinite \cite{X1,X3}, it is then enough 
to apply Theorem \ref{str.add.}.
\endproof
It then follows from \cite{KLM} that for the above coset nets the 
tensor category of all represenations is rational and modular, as shown 
in \cite{X3}, and the 
results in \cite{KLM} apply.
\subsubsection{Proof of Theorem \ref{str.add.}}
\label{proof}
The remaining and more difficult part to prove in Theorem \ref{str.add.} is
that $\A$ split and strongly additive 
implies that $\B$ is strongly additive\footnote{The reader should be 
warned that several proofs of strong additivity for specific models 
in the literature are fallacious.}. In the following we thus 
assume that $\A$ is split and strongly additive and prove that $\B$ is 
strongly additive in a series of Lemmata. The starting argument is 
similar to the one in \cite{X}.
\begin{lemma}
Let $\A$ be an irreducible, split and strongly additive, local 
conformal net on $S^1$. If $\B\subset \A$ a conformal subnet with 
$[\A:\B]<\infty$ and $I_1$ and $I_2$ are adjacent intervals, $I=\bar 
I_1\cup\bar I_2$, then $\A(I_1)\vee\B(I_2)\subset\A(I)$ is a 
finite-index irreducible inclusion of factors.
\end{lemma}
\proof
First notice that 
\[
(\A(I_1)\vee\B(I_2))'\cap\A(I)=\A(I_1)'\cap\A(I)\cap\B(I_2)'=\A(I_2)\cap\B(I_2)'=\mathbb C\ ,
\]
where $\A(I_1)'\cap\A(I)=\A(I_2)$ because $\A$ is strongly additive. 
Thus $\A(I_1)\vee\B(I_2)\subset\A(I)$ is an irreducible inclusion of 
factors. 

To show that $[\A(I):\A(I_1)\vee\B(I_2)]<\infty$ we set $I_1=[{\bf 
a},{\bf b}]$, $I_2=[{\bf b},{\bf c}]$ and take intervals $I^n_2=[{\bf 
b}_n,{\bf c}]$, where ${\rm lenght}I^n_2\nearrow{\rm lenght}I_2$.  
Then, by the split property, 
$\N_n\equiv\A(I_1)\vee\B(I^n_2)\subset\M_n\equiv\A(I_1)\vee\A(I^n_2)$ 
is isomorphic to 
$\A(I_1)\otimes\B(I^n_2)\subset\A(I_1)\otimes\A(I^n_2)$ and thus 
$[\M_n:\N_n]=\l^{-1}$, where $\l=[\A:\B§]^{-1}$.  As
\[
\N_n\nearrow\A(I_1)\vee\B(I_2),\qquad 
\M_n\nearrow\A(I_1)\vee\A(I_2)=\A(I),
\]
we have $[\A(I):\A(I_1)\vee\B(I_2)]\leq\l^{-1}$ by \cite[Proposition 
3]{KLM}.  \endproof As in Xu's proof for the group case, we consider 
two adjacent intervals $I_1$, $I_2$ and set $I=\bar I_1\cup \bar I_2$.  
Then $[\A(I):\A(I_1)\vee\B(I_2)]<\infty$, and we consider an 
expectation $\mu:\A(I)\to\A(I_1)\vee\B(I_2)$.

Then $\R(I_2)\equiv\mu (\A(I_2))$ is contained in $\A(I_1)'\cap\A(I)$ 
and the latter coincides with $\A(I_2)$ because $\A$ is strongly 
additive \cite{GLW}. Hence $\R(I_2)$ is a von Neumann algebra and
$\B(I_2)\subset\R(I_2)\subset\A(I_2)$. The following lemma is 
contained in \cite{X}.
\begin{lemma}\label{sa}
If $\R(I_2)=\A(I_2)$ then $\B$ is strongly additive.
\end{lemma}
\proof
If $\R(I_2)=\A(I_2)$ then \[
\mu(\A(I))\supset\A(I_1)\vee\R(I_2)=\A(I_1)\vee\A(I_2)=\A(I)\ ,
\]
thus $\mu$ is the identity and $\A(I_1)\vee\B(I_2)=\A(I)$. We then have
\[
\B(I)=\e_I(\A(I))=\e_I(\A(I_1)\vee\B(I_2))=\e_I(\A(I_1))\vee\B(I_2)=
\B(I_1)\vee\B(I_2)\ .
\]
\endproof
\begin{lemma} Assume $\R(I_2)=\B(I_2)$. Given intervals $L_0\subset 
L$, $\bar L_0\neq \bar L$, and $\epsilon >0$ there exists a projection 
$e\in\A(L)$ such that
\begin{gather}\label{proj}
\e_{L_0}(a)e = eae, \forall a\in\A(L_0)\ ,\\
(e\Omega,\Omega)>1 - \epsilon \ .\label{proj2}
\end{gather}
\end{lemma}
\proof
As $\B(I_2)'\cap\A(I_2)=\mathbb C$, there exists a unique expectation 
of $\A(I_2)$ onto $\B(I_2)$, thus
\[
\mu\res_{\A(I_2)} = \e_{I_2}
\]
is the vacuum preserving conditional expectation.

In order to show the Lemma we can clearly assume that $L_0$ and $L$ 
have one common boundary point. As the M\" obius group acts 
transitively on the family of three different points of $S^1$, 
property (\ref{proj},\ref{proj2}) does not depend on the choice of the 
pair $L_0\subset L$.

Let $e\neq 0$ be a projection in $B(\H)$ implementing $\mu$ namely
\[
\mu(a)e = eae,\ a\in\A(I).
\]
As $\mu$ acts identically on $\A(I_1)$, we have 
$e\in\A(I_1)'=\A(I'_1)$.

Setting $L=I'_1$, $L_0=I_2$ we then have: $L\supset L_0$ are intervals with one 
common boundary point and 
there exists a non-zero projection $e\in\A(L)$, such that the property in 
formula (\ref{proj}) holds, i.e.
\begin{equation}
e\neq 0 \ \& \ \e_{L_0}(a)e = eae,\ \forall a\in\A(L_0)\ .
	\label{impl}
\end{equation}
Clearly the above property (\ref{impl}) is a fortiori true if we replace $L$ 
with a larger interval and $L_0$ with a smaller interval.

Set $L_0=[{\bf a},{\bf b}]$, $L=[{\bf a},{\bf c}]$ and choose 
sequences of points ${\bf b}_n\in({\bf a},{\bf b})$ and ${\bf 
c}_n\in({\bf c},{\bf a})$ in $S^1$, such that ${\rm lenght}[{\bf a}, 
{\bf b}_n]\searrow 0$ and ${\rm lenght}[{\bf a},{\bf c}_n]\nearrow 
2\pi$.

As $[{\bf b}_n, {\bf c}_n]$ is an increasing sequence of intervals and 
$\cup_n [{\bf b}_n, {\bf c}_n]$ is dense in $S^1$, it follows that 
$\cup\A([{\bf b}_n, {\bf c}_n])$ is strongly dense in $B(\H)$ (this is 
a consequence of Haag duality and the factoriality of the local von 
Neumann algebras). Therefore the unitaries of $\cup_n \A([{\bf b}_n, 
{\bf c}_n])$ are strongly dense in the unitaries of $B(\H)$. Given 
$\epsilon >0$, then there exists an integer $n$ and a unitary $u\in\A([{\bf 
b}_n, {\bf c}_n])$ such that
\[
(eu\Om,u\Om)>1-\epsilon\ .
\]
Replacing $L$ with $[{\bf a}, {\bf c}_n]$, $L_0$ with $[{\bf a}, 
{\bf b}_n]$ and $e$ with $e'\equiv u^* eu$, equation 
(\ref{proj2}) clearly holds.
But also equation (\ref{proj}) is satisfied because
\[
\e_{L_0}(a)e' =u^*\e_{L_0}(a)eu= u^* eaeu=u^* euau^*eu=e'ae',
\quad a\in\A(L_0)\ ,
\]
as $u$ commutes with $\A(L_0)$.
\endproof 
\begin{lemma}\label{neq} $\R(I_2)\neq \B(I_2)$ unless $\B=\A$.
\end{lemma}
\proof 
Let's assume $\R(I_2)=\B(I_2)$. Note that by M\" obius covariance the 
equality $\R(I_2)=\B(I_2)$ is independent of the choice of $I_1,I_2$.

Let $I$ be a fixed interval and $I_n$ a decreasing sequence of 
intervals with a common boundary point with $I$ such that $\cap_n I_n 
=I$ and choose a projection $e_n\in\A(I_n)$ such that 
\begin{equation}\label{nn} \e_I(a)e_n = e_n ae_n, \forall a\in\A(I)\; 
\&\; (e_n\Om,\Om)>1-\frac{1}{n}\ .
\end{equation}
Let $e$ be a weak limit point of $\{e_n\}$. Then 
$e\in\cap\A(I_n)=\A(I)$ and $e\in\B(I)'$, thus $e$ is a 
scalar, $0\leq e\leq 1$. As $(e_n\Om,\Om)>1-\frac{1}{n}$, we 
have $e = 1$. Thus $e_n\to 1$ weakly. As the weak and strong 
topologies coincide on the set of selfadjoint projections, $e_n\to 1$ 
strongly. Going to the limit in eq. (\ref{nn}) we then have 
$\e_I(a)=a$, $a\in\A(I)$, namely $\B(I)=\A(I)$.
\endproof
\noindent
The assumptions in the following lemma will later be proven to be 
impossible.
\begin{lemma}\label{lemmass} 
Let $\A$ be a local irreducible conformal net and $\B$ a finite-index 
subnet.  Suppose $\A$ is split and strongly additive and $\B$ is not 
strongly additive.  Then there exists an intermediate conformal net 
$\B\subset\L\subset\A$ such that $\L$ is split and strongly additive 
and $\L\neq\A$.
\end{lemma}
\proof We use the above notations.  By Lemma \ref{neq} 
$\R(I_2)\neq\B(I_2)$ and Lemma \ref{sa} $\R(I_2)\neq\A(I_2)$.  By 
Theorem \ref{int-net} there exists a conformal subnet $\R$ intermediate 
between $\B$ and $\A$ such that the associated local von Neumann 
algebra $\R(I_2)$ is such a factor.  Set 
$\R_1\equiv\R$.  Replacing $\B$ by $\R_1$ and repeating the 
construction we find a factor $\R_2(I_2)$ between $\R_1(I_2)$ and 
$\A(I_2)$.  Iterating the procedure we get a sequence of 
factors $\R_n(I_2)$, coming from a conformal subnets $\R_n$, such 
that
\[
\B(I_2)\subset \R_1(I_2)\subset 
\R_2(I_2)\subset\cdots\subset \A(I_2)\ .
\]
As $[\A(I_2):\B(I_2)]<\infty$, after finitely many steps the iteration 
stabilizes, so let $n$ be the smallest integer such that $\R_n(I_2)=
\R_{n+1}(I_2)$. Then $n\geq 2$ by the above discussion. 

By Lemma \ref{neq} we then have $\R_n(I_2)= \A(I_2)$.  Thus 
$\L\equiv\R_{n-1}$ is strongly additive by Lemma \ref{sa}. $\L$ is
clearly split and, by construction, properly between $\B$ and $\A$.  
\endproof
\begin{proposition}\label{ss} 
Let $\A$ be a local irreducible conformal net 
and $\B$ a finite-index subnet. If $\A$ is split and strongly 
additive then $\B$ is split and strongly additive.
\end{proposition}
\proof 
Clearly $\B$ is split. Suppose that $\B$ is not strongly additive.
Then by Lemma 
\ref{lemmass} there exists a strongly additive
conformal subnet $\L_1$ intermediate 
between $\B$ and $\A$ such that $\L_1\neq\A$. Again  by Lemma 
\ref{lemmass} there exists a strongly additive
conformal subnet $\L_2$ intermediate 
between $\B$ and $\L_2$ such that $\L_2\neq\L_1$. By iteration
we find a sequence of subnets
\[
\A\supset\L_1\supset\L_2\supset\cdots\supset\B\ ,
\]
where all inclusions are proper, thus $[\L_{n+1} :\L_n]\geq 2$ by 
Jones' theorem \cite{J}. So $[\A:\B]=\infty$ by the multiplicativity 
of the index \cite{KL}, contrary to our assumptions.
\endproof 
{\bf Proof of Theorem \ref{str.add.}}\, It is now sufficient to apply
Lemma \ref{usplit}, Lemma \ref{uss} and Proposition \ref{ss}.
\endproof
\subsubsection{Case of a Fermi net}\label{Fermi}
Most of our analysis extends to the case of non-local Fermi conformal 
nets. As there are several examples of local conformal finite-index 
subnets of Fermi nets, we sketch here how to modify our arguments and
reduce to the local situation.

In this subsection $\A$ is a {\it twisted local} irreducible net of 
von Neumann algebras on $S^1$, namely $\A$ an irreducible net 
satisfying all properties ${\bf A}$ to ${\bf E}$ in Section \ref{nets}, 
except ${\bf B}$ which is now replaced by
\begin{itemize}
\item[${\bf B}'.$] {\it Twisted locality}.
There exists a unitary $Z$ commuting with the unitary representation 
$U$ such that $Z\Om=\Om$ and 
$$ Z\A(I')Z^*\subset \A(I)' $$
for all intervals $I$.
\end{itemize}
The basic results for local nets (modular structure, duality, etc.)
have a version for twisted local nets, see \cite{DLR}.

We shall say that a conformal net $\A$ is a {\it Fermi net} if 
there exists a self-adjoint unitary $V$ on $\H$ such that $V\Om=\Om$ 
and $\b(\A(I))=\A(I)$, $I\in\I$, where $\b\equiv{\rm Ad}V$, with 
canonical commutation relations: if $I_1,I_2$ are disjoint intervals 
then the commutator or anti-commutator
\[
[a_1,a_2]_{\pm}=0, \quad a_i\in\A(I_i)\ ,
\]
if $\b(a_i)=\pm a_i$; the commutator vanishes if one of the $a_i$'s 
is a Bose operator ($\b(a_i)=a_i$) and the anti-commutator vanishes if 
both the $a_i$ are Fermi operators ($\b(a_i)= -a_i$).

A conformal Fermi net satisfies twisted locality, hence twisted duality, 
where the unitary $Z$ is given by
\[
Z=\frac{1+iV}{1+i}\ ,
\]
in particular $ZbZ^*=b$, for all $b\in\vee_I\B_b(I)$, see \cite{GL}, 
where we denote by $\B_b$ the Bose subnet $\A^\b$ of $\A$. Note that 
$[\A:\B_b]=2$ and $\B_b$ is a local conformal net.
\begin{lemma}
Let $\B$ be a local subnet of $\A$. Then $\B\subset\B_b$.
\end{lemma}
\proof This is obvious since otherwise each $\B(I)$ would contain operators
with non-zero Fermi part and these do not commute if they are 
localized in disjoint intervals.
\endproof
Due to the above lemma, the results in the previous sections
extend to the case of a local 
finite-index subnet $\B$ of $\A$ once we show them in the particular 
case $\B=\B_b$. We give here explicitly the extension of Theorem 
\ref{str.add.}.
\begin{proposition}
Let $\A$ be an irreducible Fermi conformal net on $S^1$ and $\B\subset 
\A$ a local conformal subnet with $[\A:\B]<\infty$.
Then $\A$ is split and strongly additive iff $\B$ is split and 
strongly additive.
\end{proposition}
\proof
By the above discussion we may assume that $\B$ is the Bose subnet. 
We assume that $\A$ is split and strong additive and show that $\B$ 
strong additive, the other implications are obtained essentially as in 
the local case.

It is enough to show that $\A(I_1)\vee\B(I_2)=\A(I)$ if $I_1,I_2$ are 
adjacent intervals and $I=\bar I_1\cup \bar I_2$. The inclusion 
$\A(I_1)\vee\B(I_2)\subset\A(I)$ has finite index as in the local 
case and we assume $\A(I_1)\vee\B(I_2)\neq\A(I)$. We consider an expectation $\mu :\A(I)\to\A(I_1)\vee\B(I_2)$. 
Let $u_i\in\A(I_i)$ be Fermi unitaries. Then ${\rm Ad}u_2$ implements 
an automorphism of $\B(I_2)$, acts trivially on 
$\B(I_1)$ and ${\rm Ad}u_2(u_1)=-u_1$, therefore ${\rm Ad}u_2$ implements 
an automorphism of $\A(I_1)\vee\B(I_2)$. As $\A(I)$ is generated by
$\A(I_1)\vee\B(I_2)$ and $u_2$, it follows that 
$\A(I)$ is the crossed product of
$\A(I_1)\vee\B(I_2)$ by ${\rm Ad}u_2$, thus $\mu$ acts trivially on 
$\A(I_1)$ and $\mu(u_2)=0$, so $\mu(\A(I_2))=\B(I_2)$.

Let $e\neq 0$ be a projection implementing $\mu$. Then $e\in 
Z\A(I'_1)Z^*$. Arguing as in the local case $e\in\A(I_2)'\cap 
Z\A(I_2)Z^*$, but
\begin{multline*}
\A(I_2)'\cap Z\A(I_2)Z^* =Z(Z^*\A(I_2)'Z\cap \A(I_2))Z^*\\
\subset Z(Z^*\B(I_2)'Z\cap \A(I_2))Z^*
=Z(\B(I_2)'\cap \A(I_2))Z^*=\mathbb C \ .
\end{multline*}
so $e=1$ and $\mu$ is trivial, which contradicts the 
assumption $\A(I_1)\vee\B(I_2)\neq\A(I)$.
\endproof

\noindent {\bf Remark.} In this paper the positivity of the energy has 
been used only indirectly, essentially to entail the Reeh-Schlieder 
property and the factoriality of the local algebras.  Thus our results 
extend to the case of conformal nets on $S^1$ with the above 
properties, without assuming the positivity of the conformal 
Hamiltonian.  We encounter this situation if we consider a local 
conformal net $\A$ on the $1+1$-dimensional Minkowski spacetime and 
look at the corresponding time zero net $\A_{t=0}$.  In particular, 
if $\A$ is a finite-index local extension of a $1+1$-dimensional 
chiral net $\A_1\otimes\A_2$, as discussed in \cite{R}, then 
$\A_{t=0}$ is split and strongly additive iff both $\A_1$ and $\A_2$ 
are split and strongly additive.
\bigskip

\noindent {\bf Acknowledgements.} This work began while the author was 
visiting the Mathematical Sciences Research Institute, Berkeley, in 
November-December 2000, during the Operator Algebra program and he 
wishes to thank the organizers for the kind invitation and the warm 
hospitality.  He also thanks M. M\"{u}ger, V. Toledano and F.~Xu for 
stimulating conversations and S. Carpi and R. Conti for critical 
reading of the manuscript.

\end{document}